\magnification=1200

\input amstex

\documentstyle{amsppt}


\pagewidth{165truemm}
\pageheight{227truemm}


\def\p#1{{{\Bbb P}^{#1}_{k}}}

\def\a#1{{{\Bbb A}^{#1}_{k}}}

\def\Hilb{{{\Cal H}\kern -0.25ex{\italic ilb\/}}}

\def\Scand{{{\Cal S}\kern -0.25ex{\italic cand\/}}}

\def\Hom{{{\Cal H}\kern -0.25ex{\italic om\/}}}

\def\Ext{{{\Cal E}\kern -0.25ex{\italic xt\/}}}

\def\Sim{{{\Cal S}\kern -0.25ex{\italic ym\/}}}

\def\Ker{{{\Cal K}\kern -0.25ex{\italic er\/}}}

\def\Sing{\operatorname{Sing}}

\def\emdim{\operatorname{emdim}}

\def\gr{\operatorname{gr}}

\def\spec{\operatorname{spec}}

\def\ext{\operatorname{Ext}}

\def\lev{\operatorname{msdeg}}

\def\Soc{\operatorname{Soc}}

\def\Ofa#1{{{\Cal O}_{#1}}}

\def\M{\operatorname{\frak M}}


\topmatter

\title
Irreducibility of the Gorenstein locus of the punctual Hilbert scheme of degree $10$
\endtitle

\rightheadtext{Irreducibility of the Gorenstein locus}

\author
Gianfranco Casnati, Roberto Notari
\endauthor

\address
Gianfranco Casnati, Dipartimento di Matematica, Politecnico di Torino,
c.so Duca degli Abruzzi 24, 10129 Torino, Italy
\endaddress

\email
casnati\@calvino.polito.it
\endemail

\address
Roberto Notari, Dipartimento di Matematica \lq\lq Francesco Brioschi\rq\rq, Politecnico di Milano,
via Bonardi 9, 20133 Milano, Italy
\endaddress

\email
roberto.notari\@polimi.it
\endemail



\abstract
Let $k$ be an algebraically closed field of characteristic $0$ and let $\Hilb_{d}^{G}(\p{N})$ be the open locus of the Hilbert scheme $\Hilb_{d}(\p{N})$ corresponding to Gorenstein subschemes. We proved in a previous paper that $\Hilb_{d}^{G}(\p{N})$ is irreducible for $d\le9$ and $N\ge1$. In the present paper we prove that also $\Hilb_{10}^{G}(\p{N})$ is irreducible for each $N\ge1$, giving also  a complete description of its singular locus.
\endabstract

\endtopmatter

\document

\head
1. Introduction and notation
\endhead

Let $k$ be an algebraically closed field of characteristic $0$ and
denote by $\Hilb_{d}(\p N)$ the Hilbert scheme parametrizing
closed subschemes in $\p N$ of dimension $0$ and degree $d$.

On one hand it is well--known that  such a scheme is always
connected (see [Ha1]) and it is actually irreducible when either
$d\ge1$ and $N\le2$ (see [Fo] where a more general result is
proven) or $d\le7$ and $N\ge1$ (see [C--E--V--V]).

On the other hand, in [Ia1] the author proved that, if $d$ is
large with respect to $N$, $\Hilb_{d}(\p N)$ is always reducible.
Indeed for every $d$ and $N$ there always exists a generically
smooth component of $\Hilb_{d}(\p N)$ having dimension $dN$, the
general point of which corresponds to a reduced set of $d$ points
but, for $d$ large with respect to $N>2$, there is at least one
other component with general point corresponding to an irreducible
scheme of degree $d$ supported on a single point. For example in
the above quoted paper [C--E--V--V], the authors also prove the
existence of exactly two components in  $\Hilb_{8}(\p{N})$,
$N\ge4$.

In view of these results it is reasonable to consider the
irreducibility of other naturally occurring loci in $\Hilb_{d}(\p
N)$. E.g. one of the loci that has interested us is the set
$\Hilb_{d}^G(\p N)$ of points in $\Hilb_{d}(\p N)$ representing
schemes which are Gorenstein. This is an important locus since it
includes reduced schemes.

A first result, part of the folklore, gives the irreducibility and
smoothness of $\Hilb_{d}^G(\p N)$ when $N\le 3$. In [C--N2] (see
also [C--N1]) we proved the irreducibility of $\Hilb_{d}^G(\p N)$
when $d\le9$ and $N\ge1$. In [I--E] the authors stated that
$\Hilb_{10}^G(\p N)$ is reducible, essentially by producing an
irreducible scheme of dimension $0$ and degree $10$ corresponding
to a point in the Hilbert scheme having tangent space of too small
dimension. Unfortunately their computations where affected by a
numerical mistake as R. Buchweitz pointed out.  In [I--K], Lemma
6.21, the authors claim the reducibility of $\Hilb_{14}^G(\p 6)$,
asserting the existence of numerical examples that can be checked
using the \lq\lq Macaulay\rq\rq algebra program.

The main result of this paper is the following

\proclaim{Main Theorem} Let $k$ be an algebraically closed field
of characteristic $0$. Then the scheme $\Hilb_{10}^{G}(\p{N})$  is
irreducible for each $N\ge1$. \qed
\endproclaim

In order to prove the above theorem we will also make use of the
classification results proved in [C--N2] and [Cs]. The proof of
the Main Theorem is given in Section 4. It rests on the analysis
of several different cases, which we examine separately in
Sections 2, 3, 4.

The idea is that each $X\in \Hilb_{10}^{G}(\p{N})$ is the spectrum
of an Artinian Gorenstein $k$--algebra $A$ and the irreducibility
of $\Hilb_{10}^{G}(\p{N})$ depends on some properties of $A$ which
can be checked on the direct summands of $A$, which correspond to
the irreducible components of the original scheme $X$. Thus we can
restrict our attention to local algebras $A$ with maximal ideal
$\frak M$, using all the known classification results.

More precisely in Section 2 we list some preliminary results. In
particular we recall that the algebras which we are interested in
satisfy $\dim_k(\M/\M^2)\le4$. In Section 3 we examine Artinian,
Gorenstein local $k$--algebras of degree $d\le10$ for which
$\dim_k(\M^2/\M^3)\le3$,  with the same methods used in [C--N1],
[C--N2] and [Cs]. Artinian, Gorenstein local $k$--algebras of
degree $d\le10$ with $\dim_k(\M^2/\M^3)=4$ cannot be easily
treated in this way so, in Section 4, we analyse this remaining
case via an indirect approach.

In the last Section 5 we deal with the singular locus of
$\Hilb_{10}^{G}(\p{N})$. Again the fact that
$X\in\Hilb_{10}^{G}(\p{N})$ is singular  in its Hilbert scheme (we
briefly say that $X$ is obstructed in this case) can be recovered
from the local direct summands of the associated algebra $A$. As
for the irreducibility we are able to give an easy criterion for
deciding weather a fixed scheme  $X$ is obstructed or not in term
of the underlying algebra.

\subhead
Notation
\endsubhead
In what follows $k$ is an algebraically closed field of
characteristic $0$.

Recall that a Cohen--Macaulay local ring $R$ is one for which
$\dim(R)={\roman{depth}}(R)$. If, in addition, the injective
dimension of $R$ is finite then $R$ is called Gorenstein
(equivalently, if $\ext_R^i\big(M,R)=0$ for each $R$--module $M$
and $i>\dim(R)$).  An arbitrary ring $R$ is called Cohen--Macaulay
(resp. Gorenstein) if $R_{\frak M}$ is Cohen--Macaulay (resp.
Gorenstein) for every maximal ideal ${\frak M}\subseteq R$.

All the schemes $X$ are separated and of finite type over $k$. A
scheme $X$ is Cohen--Macaulay (resp. Gorenstein) if for each point
$x\in X$ the ring ${\Cal O}_{X,x}$ is Cohen--Macaulay (resp.
Gorenstein). The scheme $X$ is Gorenstein if and only if it is
Cohen--Macaulay and its dualizing sheaf $\omega_{X}$ is
invertible.

For each numerical polynomial $p(t)\in{\Bbb Q}[t]$, we denote by
$\Hilb_{p(t)}(\p N)$ the Hilbert scheme of closed subschemes of
$\p N$ with Hilbert polynomial $p(t)$. With abuse of notation we
will denote by the same symbol both a point in $\Hilb_{p(t)}(\p
N)$ and the corresponding subscheme of $\p N$. In particular we
will say that $X$ is obstructed (resp. unobstructed) in $\p N$ if
the corresponding point is singular (resp. non--singular) in
$\Hilb_{p(t)}(\p N)$.

Moreover we denote by $\Hilb_{p(t)}^G(\p N)$ the locus of points
representing Gorenstein schemes. This is an open subset of
$\Hilb_{p(t)}(\p N)$, though not necessarily dense.

If $X\subseteq\p N$ we will denote by $\Im_X$ its sheaf of ideals
in $\Ofa{\p N}$ and we define the normal sheaf of $X$ in $\p N$ as
${\Cal N}_X:=(\Im_X/\Im_X^2)\check{\
}:=\Hom_{X}\big(\Im_X/\Im_X^2,\Ofa X\big)$. If we wish to stress
the fixed embedding $X\subseteq\p N$ we will write ${\Cal
N}_{X\vert\p N}$ insted of ${\Cal N}_X$. If $X\in\Hilb_{p(t)}(\p
N)$, the space $H^0\big(\p N,{\Cal N}_X\big)$ can be canonically
identified with the tangent space to $\Hilb_{p(t)}(\p N)$ at the
point $X$. In particular $X$ is obstructed in $\p N$ if and only
if $h^0\big(\p N,{\Cal N}_X\big)$ is greater than the local
dimension of $\Hilb_{p(t)}(\p N)$ at the point $X$.

If $\gamma:=(\gamma_0,\dots,\gamma_n)\in{\Bbb N}^{n+1}$ is a
multi--index, then we set $\vert
\gamma\vert:=\sum_{i=0}^n\gamma_i$,
$\gamma!:=\prod_{i=0}^n{\gamma_i!}$,
$t^\gamma:=t_0^{\gamma_0}\dots t_n^{\gamma_n}\in k[t_0,\dots,t_n]$
and we say that $\gamma\ge0$ if and only if $\gamma_i\ge0$ for
each $i=0,\dots,n$. If $\delta:=(\delta_0,\dots,\delta_n)\in{\Bbb
N}^{n+1}$ is another multi--index then we write $\gamma\ge\delta$
if and only if $\gamma-\delta\ge0$. Finally we set
$$
{\gamma\choose\delta}:={{\gamma!}\over{\delta!(\gamma-\delta)!}}.
$$

For all the other notations and results we refer to [Ha2].

\head
2.  Reduction to the local case
\endhead
We begin this section by recalling some general facts about
$\Hilb_{d}(\p{N})$. The locus of reduced schemes ${\Cal
R}\subseteq\Hilb_{d}(\p{N})$ is birational to a suitable open
subset of the $d$--th symmetric product of $\p N$, thus it is
irreducible of dimension $dN$ (see [Ia1]). We will denote by
$\Hilb_{d}^{gen}(\p{N})$ its closure in $\Hilb_{d}(\p{N})$

Notice that $\Hilb_{d}^{gen}(\p{N})$ is necessarily an irreducible
component of $\Hilb_{d}(\p{N})$. Indeed, in any case, we can
always assume $\Hilb_{d}^{gen}(\p{N})\subseteq{\Cal H}$ for a
suitable irreducible component $\Cal H$ in $\Hilb_{d}(\p{N})$. If
the inclusion were proper then there would exist a flat family
with special point in ${\Cal R}$, hence reduced, and non--reduced
general point, which is absurd. We conclude that
$\Hilb_{d}^{gen}(\p{N})={\Cal H}$.

\definition{Definition 2.1}
A scheme $X$ is said to be smoothable in $\p N$ if
$X\in\Hilb_{d}^{gen}(\p{N})$.
\enddefinition

Thus $X$ is smoothable if and only if there exists an irreducible
scheme $B$ and a flat family ${\Cal X}\subseteq\p N\times B\to B$
with special fibre $X$ and general fibre in $\Cal R$, hence
reduced. Moreover it is clear that $X$ is smoothable if and only
if the same is true for all its connected components (which
coincide with its  irreducible components since $X$ has dimension
$0$).

The following result is well--known (see e.g. [C--N2], Lemma 2.2).

\proclaim{Lemma 2.2} Let $X$ be a scheme of dimension $0$ and
degree $d$ and let $X\subseteq\p N$ and $X\subseteq\p{N'}$ be two
embeddings. Then $X$ is smoothable in $\p N$ if and only if it is
smoothable in $\p{N'}$. \qed
\endproclaim

We now quickly turn our attention to the singular locus of
$\Hilb_{d}(\p{N})$. We have (see e.g. [C--N2], Lemma 2.3)

\proclaim{Lemma 2.3} Let $X$ be a scheme of dimension $0$ and
degree $d$ and let $X\subseteq\p N$ and $X\subseteq\p{N'}$ be two
embeddings. Then
$$
h^0\big(X,{\Cal N}_{X\vert\p N}\big)-dN=h^0\big(X,{\Cal N}_{X\vert
\p {N'}}\big)-dN'.\qquad\qed
$$
\endproclaim

Thanks to the Lemma above it follows that the obstructedness of
$X\in\Hilb_{d}^{gen}(\p{N})$ can be checked with respect to an
arbitrary embedding $X\subseteq\p N$. Moreover, if
$X=\bigcup_{i=1}^pX_i$ where $X_i$ is irreducible of degree $d_i$,
then
$$
h^0\big(\p N,{\Cal N}_X\big)=\sum_{i=1}^ph^0\big(\p N,{\Cal N}_{X_i}\big),\tag 2.4
$$
thus $X$ is unobstructed if and only if the same is true for all
its components $X_i$.

Now we restrict to $X\in\Hilb_{d}^{G}(\p{N})\subseteq
\Hilb_{d}(\p{N})$ the Gorenstein locus, i.e. the locus of points
in $\Hilb_{d}(\p{N})$ representing Gorenstein subschemes of $\p
N$. Such a locus  is actually open inside $\Hilb_{d}(\p N)$, since
its complement coincides with the locus of points over which the
relative dualizing sheaf of the universal family is not
invertible. However the locus $\Hilb_{d}^{G}(\p{N})$ is not
necessarily dense.

Trivially ${\Cal R}\subseteq\Hilb_{d}^{G}(\p{N})$, i.e. reduced
schemes represent points in $\Hilb_{d}^{G}(\p{N})$. It follows
that the main component
$\Hilb_{d}^{G,gen}(\p{N}):=\Hilb_{d}^G(\p{N})\cap\Hilb_{d}^{gen}(\p{N})$
of $\Hilb_{d}^G(\p{N})$ is irreducible of dimension $dN$ and open
in $\Hilb_{d}^G(\p{N})$ since $\Hilb_{d}^G(\p{N})$ is open in
$\Hilb_{d}(\p{N})$ (see the introduction).

As first step in the description of $\Hilb_{d}^G(\p{N})$ we show
that we can restrict our attention to schemes
$X\in\Hilb_{d}^{G}(\p{N})$ having \lq\lq big\rq\rq\  tangent space
at some point. More precisely we have the following (see e.g.
[C--N2], Proposition 2.5).

\proclaim{Proposition 2.5} Let $X\in\Hilb_d^G(\p N)$. If the
dimension of the tangent space at every point of $X$ is at most
three, then $X\in\Hilb_d^{G,gen}(\p N)$ and it is unobstructed.
\endproclaim

In order to study the obstructedness of $X\in\Hilb_d^G(\p N)$ we
finally recall that
$$
h^0\big(X,{\Cal N}_{X}\big)=\deg(X^{(2)})-\deg(X)\tag2.6
$$
where $X^{(2)}$ is the first infinitesimal neighborhood of $X$ in
$\p N$ (see Proposition 5.5 of [C--N2]).

From now on we turn our attention from $d$ general to $d=10$, i.e.
we consider  $\Hilb_{10}^{G}(\p{N})$. In order to prove its
irreducibility it thus suffices to prove the equality
$\Hilb_{10}^{G}(\p{N})=\Hilb_{10}^{G,gen}(\p{N})$, i.e. that each
$X\in \Hilb_{10}^{G}(\p{N})$ is smoothable.

Since we proved in [C--N2] that $\Hilb_d^G(\p N)$ is irreducible
if $d\le9$ and smoothability can be checked componentwise, we
deduce the following

\proclaim{Proposition 2.7} Let $X\in\Hilb_d^G(\p N)$. If all the
irreducible components of $X$ have degree at most $9$, then
$X\in\Hilb_d^{G,gen}(\p N)$. \qed
\endproclaim

In order to complete the proof of the Main Theorem stated in the
introduction, thanks to Propositions 2.5 and 2.7 we thus have to
restrict our attention to irreducible schemes $X$ of degree $d=
10$ with tangent space of dimension $n\ge4$

Each such scheme is isomorphic to $\spec(A)$, where $A$ is a
suitable local, Artinian, Gorenstein $k$--algebra of degree $d=
10$ and $\emdim(A)=n\ge4$. Thus we will first recall some results
about such kind of objects.

Let $A$ be a local, Artinian $k$--algebra of degree $d$ with
maximal ideal $\M$. In general we have a filtration
$$
A\supset\M\supset\M^2\supset\dots\supset\M^e\supset\M^{e+1}=0
$$
for some integer $e\ge1$, so that its associated graded algebra
$$
\gr(A):=\bigoplus_{i=0}^\infty\M^i/\M^{i+1}
$$
is a vector space over $k\cong A/\M$ of finite dimension
$d=\dim_k(A)=\dim_k(\gr(A))=\sum_{i=0}^e\dim_k(\M^i/\M^{i+1})$.
The Hilbert function of $A$ is by definition the function
$h_A\colon{\Bbb N}\to{\Bbb N}$ defined by
$h_A(i):=\dim_k(\M^i/\M^{i+1})$.

We recall the definition of the {\sl maximum socle degree}\/ of a
local, Artinian $k$--algebra.

\definition{Definition 2.8}
Let $A$ be a local, Artinian $k$--algebra. If $\M^e\ne0$ and
$\M^{e+1}=0$ we define the maximum socle degree of $A$ as $e$ and denote it by
$\lev(A)$.

If $e=\lev(A)$ and $n_i:=\dim_k(\M^i/\M^{i+1})$, $0\le i\le e$,
then the Hilbert function $h_A$ of $A$ will be often identified
with the vector $(n_0,\dots,n_e)\in{\Bbb N}^{e+1}$.
\enddefinition

In any case $n_0=1$. Recall that the Gorenstein condition is
equivalent to saying that the socle $\Soc(A):=0\colon \M$ of $A$
is a vector space over $k\cong A/\M$ of dimension $1$. If
$e=\lev(A)\ge1$ trivially $\M^e\subseteq\Soc(A)$, hence if $A$ is
Gorenstein then equality must hold and $n_e=1$, thus if
$\emdim(A)\ge2$ we deduce that $\lev(A)\ge2$ and
$\deg(A)\ge\emdim(A)+2$.

Taking into account of Section 5F of [Ia4] (see also [Ia2]), the
list of all possible shapes of Hilbert functions of local,
Artinian, Gorenstein $k$--algebra $A$ of degree $d=10$ and
$\emdim(A)\ge4$ is
$$
\gathered
(1,4,1,1,1,1,1),\ (1,5,1,1,1,1),\ (1,6,1,1,1),\ (1,7,1,1),\ (1,8,1)\\
(1,4,2,1,1,1),\ (1,4,2,2,1),\ (1,5,2,1,1),\ (1,6,2,1)\\
(1,4,3,1,1),\ (1,5,3,1),\\
(1,4,4,1).
\endgathered\tag 2.9
$$
As we will see later on all the above sequences actually occur as
Hilbert functions of some local, Artinian, Gorenstein
$k$--algebra. They can be divided into four different families
according to $\dim_k(\M^2/\M^3)$.

In the next two sections we will examine separately the two cases
$\dim_k(\M^2/\M^3)\le3$ and $\dim_k(\M^2/\M^3)=4$, completing the
proof of the Main Theorem.

\head
3. The cases $\dim_k(\M^2/\M^3)\le3$
\endhead

When $\dim_k(\M^2/\M^3)=1$ the sequences on the first line of
(2.9) completely characterize the algebra (see [Sa]; another proof
can be found in [C--N]), since for a local, Artinian $k$--algebra
$A$ of degree $d\ge n+2$, one has $h_A=(1,n,1,\dots,1)$ if and
only if $A\cong A_{n,d}$ where
$$
A_{n,d}:=\cases
k[x_1]/(x_1^{d})&\text{if $n=1$},\\
k[x_1,\dots,x_n]/(x_ix_j, x_h^2-x_1^{d-n})_{1\le i<j\le n,\atop
2\le h\le n}&\text{if $n\ge2$}.
\endcases
$$
Moreover we have

\proclaim{Proposition 3.1} Let $X\cong\spec(A_{n,d})\subseteq\p
N$, $N\ge n$. Then $X$ is smoothable in $\p N$.
\endproclaim
\demo{Proof} By induction on $d$, it suffices to show that
$A_{n,d}$ is a flat specialization of the simpler algebra
$A_{n,d-1}\oplus A_{0,1}$, for each $d\ge n+2\ge4$ and we refer
the reader to Remark 2.10 of [C--N2] for the details. \qed
\enddemo

We now go to examine the case $\dim_k(\M^2/\M^3)=2$, i.e. we are
considering the sequences on the second line of (2.9). If
$h_A=(1,n,2,1,\dots,1)$  (hence $\dim_k(\M^3/\M^4)=1$) it has been
already described in [E--V] (see also Section 3 of [C--N2]). In
particular $A\cong A_{n,2,d}^t:=k[x_1,\dots,x_n]/I_t$, $t=1,2$,
where
$$
\gather
{I_1:=\cases
(x_1^2x_2-x_1^3,x_2^2,x_ix_j,x_h^2-x_1^{3})_{{1\le i<j\le n,\ 3\le j}\atop 3\le h\le n} &\text{if $d=n+4$,}\\
(x_1^2x_2,x_2^2-x_1^{d-n-2},x_ix_j,x_h^2-x_1^{d-n-1})_{{1\le i<j\le n,\ 3\le j}\atop 3\le h\le n}&\text{if $d\ge n+5$,}
\endcases}\\
I_2:=(x_1x_2,x_2^3-x_1^{d-n-1},x_ix_j,x_h^2-x_1^{d-n-1})_{{1\le i< j\le n,\ 3\le j}\atop 3\le h\le n}.
\endgather
$$
Also in this case we have

\proclaim{Proposition 3.2} Let
$X\cong\spec(A_{n,2,d}^t)\subseteq\p N$, $N\ge n$. Then $X$ is
smoothable in $\p N$.
\endproclaim
\demo{Proof}
See Remark 3.4 of [C--N2].
\qed
\enddemo

If $h_A=(1,4,2,2,1)$ (hence $\dim_k(\M^3/\M^4)=2$) the algebra $A$
can be easily described making use of [Cs], Section 4. In this
case $A\cong A_{4,2,2,10}^t:=k[x_1,\dots,x_n]/I_t$, $t=1,2,3$,
where
$$
\gather
I_1:=(x_1x_2,x_2^{4}-x_1^{4},x_ix_j,x_j^2-x_1^{4})_{{1\le i<j\le 4}\atop 3\le j},\\
I_2:=(x_1^{3}x_2-x_1^{4},x_2^2,x_ix_j,x_j^2-x_1^{{4}})_{{1\le i<j\le 4}\atop 3\le j},\\
I_3:=(x_1^{3}x_2-x_1^{4},x_2^2-x_1^{3},x_ix_j,x_j^2-x_1^{{4}},x_1^5)_{{1\le i<j\le 4}\atop 3\le j}.
\endgather
$$

Also in this case we have

\proclaim{Proposition 3.3} Let
$X\cong\spec(A_{4,2,2,10}^t)\subseteq\p N$, $N\ge 4$. Then $X$ is
smoothable in $\p N$.
\endproclaim
\demo{Proof} We will give explicit flat families with general
fibre in $ \Hilb_{10}^{G,gen}(\p{N})$ and special fibre isomorphic
to $\spec(A_{4,2,2,10}^t)$, $t=1,2,3$. To this purpose take
$$
\gather
J_1:=(x_1x_2,x_2^{4}-x_1^{4},x_ix_j,x_3^2-x_1^{4},x_4^2-bx_4-x_1^{4})_{{1\le i<j\le 4}\atop 3\le j},\\
J_2:=(x_1^{3}x_2-x_1^{4},x_2^2,x_ix_j,x_3^2-x_1^{{4}},x_4^2-bx_4-x_1^{4})_{{1\le i<j\le 4}\atop 3\le j},\\
J_3:=(x_1^{3}x_2-x_1^{4},x_2^2-x_1^{3},x_ix_j,x_3^2-x_1^{{4}},x_4^2-bx_4-x_1^{4},x_1^5)_{{1\le i<j\le 4}\atop 3\le j}.
\endgather
$$
We claim that the family ${\Cal
A}^{t}:=k[b,x_1,x_2,x_3,x_4]/J_t\to \a1$ has special fibre over
$b=0$ isomorphic to  $A_{4,2,2,10}^{t}$ and general fibre
isomorphic to $A_{3,2,2,9}^{t}\oplus A_{0,1}$. In particular the
family ${\Cal A}^{t}$ is flat and has general fibre in
$\Hilb_{10}^{G,gen}(\p{N})$ due to  Proposition 2.7, thus it turns
out that also its special fibre $X$ is  in
$\Hilb_{10}^{G,gen}(\p{N})$.

It thus remains to prove the claim. To this purpose let us examine
only the case $t=1$, the other ones being similar. Let
$$
J_0:=(x_1,x_2,x_3,x_4-b)\cap(x_1x_2,x_2^{4}-x_1^{4},x_ix_j,x_3^2-x_1^{4},x_4^2, bx_4+x_1^{4})_{{1\le i<j\le 4}\atop 3\le j}.
$$
The inclusion $J_1\subseteq J_0$ is obvious. Conversely let $y\in
J_0$. Then
$$
y=u_1(x_2^4-x_1^4)+u_2(x_3^2-x_1^4)+u_3x_1x_2+
\sum_{{1\le i<j\le 4}\atop 3\le j}u_{i,j}x_ix_j
+vx_4^2+w(bx_4+x_1^{4})
$$
where $u_h,u_{i,j},vx_4^2,w\in k[b,x_1,x_2,x_3,x_4]$, $h=1,2,3$,
${1\le i<j\le 4}$ and $3\le j$, with the obvious condition
$vx_4^2+wbx_4\in(x_1,x_2,x_3,x_4-b)$. Since
$x_4\not\in(x_1,x_2,x_3,x_4-b)$ it follows that
$vx_4+wb\in(x_1,x_2,x_3,x_4-b)$. With a proper change of the
coefficients we can actually assume that $v,w\in k[b,x_4]$ whence
we finally obtain $w=-v$, i.e.
$$
y=u_1(x_2^4-x_1^4)+u_2(x_3^2-x_1^4)+u_3x_1x_2+
\sum_{{1\le i<j\le 4}\atop 3\le j}u_{i,j}x_ix_j
+v(x_4^2-bx_4-x_1^{4})
$$
i.e. $y\in J_1$.
\qed
\enddemo

Now we consider the case $\dim_k(\M^2/\M^3)=3$, i.e. we are
considering the sequences on the third line of (2.9). This has
been already described in Section 4 of [C--N2] when $d=n+5$, i.e.
$h_A=(1,n,3,1)$. In particular $A\cong
A_{n,3,n+5}^{t,\alpha}:=k[x_1,\dots,x_n]/I_{t,\alpha}$,
$t=1,\dots,6$, where
$$
\gather
I_{1,\alpha}:=(x_1x_2+x_3^2,x_1x_3,x_2^2-\alpha x_3^2+x_1^{2},x_ix_j,x_j^2-x_1^3)_{{1\le i<j\le n}\atop 4\le j},\\
I_{2,0}:=(x_1^2,x_2^2,x_3^2+2x_1x_2,x_ix_j,x_j^2-x_1x_2x_3)_{{1\le i<j\le n}\atop 4\le j},\\
I_{3,0}:=(x_1^2,x_2^2,x_3^2,x_ix_j,x_j^2-x_1x_2x_3)_{{1\le i<j\le n}\atop 4\le j},\\
I_{4,0}:=(x_2^3-x_1^3,x_3^3-x_1^3,x_ix_j,x_j^2-x_1^3)_{{1\le i<j\le n}\atop 4\le j},\\
I_{5,0}:=(x_1^2,x_1x_2,x_2x_3,x_2^3-x_3^3,x_1x_3^2-x_3^3,x_ix_j,x_j^2-x_3^3)_{{1\le i<j\le n}\atop 4\le j},\\
I_{6,0}:=(x_1^2,x_1x_2,2x_1x_3+x_2^2,x_3^3,x_2x_3^2,x_ix_j,x_j^2-x_1x_3^2)_{{1\le i<j\le n}\atop 4\le j}.
\endgather
$$
Also in this case we have

\proclaim{Proposition 3.4} Let
$X\cong\spec(A_{n,3,n+5}^{t,\alpha})\subseteq\p N$, $N\ge n$. Then
$X$ is smoothable in $\p N$.
\endproclaim
\demo{Proof}
See Remark 4.9 of [C--N2].
\qed
\enddemo

If $h_A=(1,4,3,1,1)$, the algebra $A$ can be described making use
of [Cs], Section 5. In this case $A\cong
A_{4,3,10}^t:=k[x_1,\dots,x_n]/I_t$, $t=0,\dots,6$, where
$$
\gather
I_0:=(x_1x_2+x_3^2,x_1x_3,x_2^2-x_1^{3},x_ix_4,x_4^2-x_1^4)_{{1\le i\le 3}},\\
I_1:=(x_1x_2+x_3^2,x_1x_3,x_2^2-x_3^2-x_1^{3},x_ix_4,x_4^2-x_1^4)_{{1\le i\le 3}},\\
I_2:= (x_1x_2,x_1^2-x_3^3,x_2^2-x_3^3, x_1x_3^2,x_2x_3^2, x_ix_4,x_4^2-x_3^4)_{{1\le i\le 3}},\\
I_3:=(x_1x_2,x_2x_3,x_1^2-x_3^3,x_1x_3^2,x_2^3-x_3^4,x_ix_4,x_4^2-x_3^4)_{{1\le i\le 3}},\\
I_4:=(x_1x_2,x_1x_3,x_2x_3,x_2^3-x_1^4,x_3^3-x_1^4,x_ix_4,x_4^2-x_1^4)_{{1\le i\le 3}},\\
I_5:=(x_1x_2,x_2x_3,x_1^2,x_1x_3^2-x_2^4,x_3^3-x_2^4,x_ix_4,x_4^2-x_2^4)_{{1\le i\le 3}},\\
I_6:=(x_1x_2-x_3^3,2x_1x_3+x_2^2,x_1^2,x_1x_3^2,x_2x_3^2,x_ix_4,x_4^2-x_3^4)_{{1\le i\le 3}}.
\endgather
$$

Again we have

\proclaim{Proposition 3.5} Let
$X\cong\spec(A_{4,3,10}^t)\subseteq\p N$, $N\ge 4$. Then $X$ is
smoothable in $\p N$.
\endproclaim
\demo{Proof} The argument is the same of the proof of Proposition
3.2. Indeed it suffices to take
$$
\gather
J_0:=(x_1x_2+x_3^2,x_1x_3,x_2^2-x_1^{3},x_ix_4,x_4^2-bx_4-x_1^4)_{{1\le i\le 3}},\\
J_1:=(x_1x_2+x_3^2,x_1x_3,x_2^2-x_3^2-x_1^{3},x_ix_4,x_4^2-bx_4-x_1^4)_{{1\le i\le 3}},\\
J_2:= (x_1x_2,x_1^2-x_3^3,x_2^2-x_3^3, x_1x_3^2,x_2x_3^2, x_ix_4,x_4^2-bx_4-x_3^4)_{{1\le i\le 3}},\\
J_3:=(x_1x_2,x_2x_3,x_1^2-x_3^3,x_1x_3^2,x_2^3-x_3^4,x_ix_4,x_4^2-bx_4-x_3^4)_{{1\le i\le 3}},\\
J_4:=(x_1x_2,x_1x_3,x_2x_3,x_2^3-x_1^4,x_3^3-x_1^4,x_ix_4,x_4^2-bx_4-x_1^4)_{{1\le i\le 3}},\\
J_5:=(x_1x_2,x_2x_3,x_1^2,x_1x_3^2-x_2^4,x_3^3-x_2^4,x_ix_4,x_4^2-bx_4-x_2^4)_{{1\le i\le 3}},\\
J_6:=(x_1x_2-x_3^3,2x_1x_3+x_2^2,x_1^2,x_1x_3^2,x_2x_3^2,x_ix_4,x_4^2-bx_4-x_3^4)_{{1\le i\le 3}},
\endgather
$$
observing again that ${\Cal A}^{t}:=k[b,x_1,x_2,x_3,x_4]/J_t\to
\a1$ is flat, it has special fibre over $b=0$ isomorphic to
$A_{4,3,10}^{t}$ and general fibre isomorphic to
$A_{3,3,9}^{t}\oplus A_{0,1}$. \qed
\enddemo

\remark{Remark 3.6} In Section 4 of [Cs], local, Artinian,
Gorenstein $k$--algebras $A$ with Hilbert function
$h_A=(1,n,2,\dots,2,1)$ are completely classified. Taking into
account of such a classification, it is trivial to modify the
above explicit proof of Proposition 3.3 in order to prove that
every scheme $X\cong\spec(A)$ with $h_A=(1,n,2,\dots,2,1)$ is
smoothable for each $n\ge2$.

Similarly, it is trivial to modify the proof Proposition 3.5 in
order to prove that every scheme $X\cong\spec(A)$ with
$h_A=(1,n,3,1,\dots,1)$ is smoothable for each $n\ge3$.
\endremark
\medbreak

\head
4. The case $\dim_k(\M^2/\M^3)=4$
\endhead

In this section we deal with the last case, namely
$\dim_k(\M^2/\M^3)=4$ or, equivalently, $h_A=(1,4,4,1)$. In this
case we will not exploit any explicit description for such
algebras as we did in the case $\dim_k(\M^2/\M^3)\le3$ but we will
make use of some classical results about Artinian Gorenstein
$k$--algebras combined with a recent structure Theorem for such
algebras $A$ with $h_A=(1,N,N,1)$ (see [E--R]).

Indeed on one hand Theorem 3.3 of [E--R] states that each
Artinian, Gorenstein $k$--algebras with $h_A=(1,N,N,1)$ is
canonically graded, i.e. there exists an homogeneous ideal
$I\subseteq S:=k[x_1,\dots,x_N]$ such that $A\cong S/I$.

On the other hand, in order to construct such graded quotient
algebras it suffices to make use of the theory of inverse systems
that we are going to recall very quickly (as reference see [I--K],
Section 1). We have an action of $S:=k[x_1,\dots,x_N]$ over
$R:=k[y_1,\dots,y_N]$ given by partial derivation by identifying
$x_i$ with $ \partial/ \partial {y_i}$.  Hence
$$
x^{\alpha}(y^{\beta}):=\cases
\alpha!{\beta\choose\alpha}y^{\beta-\alpha}\qquad&\text{if $\beta\ge\alpha$,}\\
0\qquad&\text{if $\beta\not\ge\alpha$.}
\endcases
$$

Such an action defines a perfect pairing $S_d\times R_d\to k$
between forms of degree $d$ in $R$ and in $S$. We will say that
{two \sl homogeneous forms $g\in R$ and $f\in S$ are apolar if
$f(g)=0$}\/. As explained in [I--K] apolarity allows us to
associate an Artinian Gorenstein graded quotient of $S$ to a form
in $R$ as follows. Let $g\in R_d$: then we set
$$
g^\perp:=\{\ f\in S\ \vert\ f(g) =0\ \}
$$
and it is easy to prove that both $g^\perp$ is a homogeneous ideal
in $S$ and $S/g^\perp$ is an Artinian Gorenstein graded quotient
of $S$ with socle in degree $d$. Also the converse is true i.e. if
$A$ is an Artinian Gorenstein graded quotient of $S$, say
$A:=S/I$, with socle in degree $d$ then there exists $g\in R_d$
such that $I=g^\perp$. The main result about apolarity due to
Macaulay (see [I--K], Lemma 2.12 and the references cited there)
is the following

\proclaim{Theorem 4.1} The map $g\mapsto S/g^\perp$ induces a
bijection between ${\Bbb P}(R_d)$ and the set of graded Artinian
Gorenstein quotient rings of $S$ with socle in degree $d$. \qed
\endproclaim

Moreover the set of polynomials corresponding to algebras $A$
having maximal embedding dimension $h_A(1)=N$ is a non--empty open
subset of ${\Bbb P}(R_d)$ due to the following standard and
well--known

\proclaim{Lemma 4.2} Let $g\in R_d$, $A:=S/g^\perp$, $t\le N$.
Then $h_A(1)\le t$ if and only if there exist
$\ell_1,\dots,\ell_t\in R_1$ such that $g\in
k[\ell_1,\dots,\ell_t]$.
\endproclaim
\demo{Proof} If $t=N$ there is nothing to prove. Assume that $
t<N$. If $h_A(1)\le t$, up to a proper change of the coordinates
$x_1,\dots,x_N\in S_1$ we can assume that $x_N\in g^\perp$, thus
$g\in k[y_1,\dots,y_{N-1}]$. Conversely if there exist
$\ell_1,\dots,\ell_t\in R_1$ such that $g\in
k[\ell_1,\dots,\ell_t]$, since $\dim_k(S_1)=N$, it follows the
existence of linear forms $\ell_{t+1},\dots,\ell_N\in S_1$ which
are not in the space spanned by $\ell_1,\dots,\ell_t$: in
particular $\ell_i(g)=0$ for such $N-t$ forms. Thus
$\ell_{t+1},\dots,\ell_N\in g^\perp$, whence
$h_A(1)=\dim_k(S_1)-\dim_k(g^\perp\cap S_1)\le t$ \qed
\enddemo

Now, we restrict our attention to algebras with Hilbert function
$(1,4,4,1)$. Thus there exists a natural variety $\Cal Z$ which
parametrizes such kind of algebras.  More precisely ${\Cal Z}$ is
the open non--empty subset of ${\Bbb P}(R_3)\cong\p{19}$ of cubic
surfaces in $\p3$ which are not cones due to the previous lemma.

From now on we will denote by ${\Cal Z}_N$ the locus of
irreducible schemes $X\in\Hilb_{10}^{G}(\p{N})$ of the form
$X=\spec(A)\subseteq\p N$ with $h_A=(1,4,4,1)$. Necessarily
$N\ge4$ and ${\Cal Z}_4={\Cal Z}$ thus it is irreducible.

Our aim is to prove that ${\Cal
Z}_N\subseteq\Hilb_{10}^{G,gen}(\p{N})$. Let us examine first the
case $N=4$. If the closure of  ${\Cal Z}_4$ in
$\Hilb_{10}^{G}(\p{4})$ were contained in a component different
from $\Hilb_{10}^{G,gen}(\p{4})$, then each smoothable $X\in{\Cal
Z}_4$, if any, would be obstructed.

In [I--E] the authors asserted the existence of such a scheme but
their computations were affected by a mistake pointed out to the
authors by R. Buchweitz in a private communication. In Example 4.1
of [Ia3] the author claimed the smoothability of all points in
${\Cal Z}_N$ without providing any proof for this. We will give
here a quick proof of this fact.

\proclaim{Proposition 4.3}
There exists an unobstructed $X\in{\Cal Z}_4\cap\Hilb_{10}^{G}(\p{4})$.
\endproclaim
\demo{Proof}
Consider the ideal
$$
J:=(x_3x_4, x_2x_4, x_1x_4, x_1^2 + x_2^2, x_1x_2 + x_3^2, x_1 x_3, x_4^3 - b^2x_4+(b-1)x_1^3, x_3^3, x_2^2x_3, x_2^3)
$$
in $k[b,x_1,x_2,x_3,x_4]$. Let ${\Cal A}:=k[b,x_1,x_2,x_3,x_4]/J$
and denote by ${\Cal X}\to \a1$ the corresponding family.

If $b\ne0$, then
$J=J_1\cap J_2$
where
$$
\gather
J_1:=(x_4^2, x_3x_4, x_2x_4, x_1x_4, x_1x_3, x_1x_2 + x_3^2, x_1^2 + x_2^2, x_3^3, x_2^2x_3, x_2^3, bx_2x_3^2 - x_2x_3^2 - b^2x_4),\\
J_2=(x_1,x_2,x_3,x_4^2-b^2)
\endgather
$$
(one can use any computer algebra system for checking such an
equality). Since we have $bx_2x_3^2 - x_2x_3^2 - b^2x_4\in J_1$,
when $b\ne0$ we have an isomorphism
$$
k[x_1,x_2,x_3,x_4]/J_1\cong k[x_1,x_2,x_3]/(x_1x_3, x_1x_2 + x_3^2, x_1^2 + x_2^2)\cong A_{3,3,8}^{1,0}.
$$
Such an algebra is smoothable by Lemma 3.4. Since the fibres
${\Cal X}_b$ with $b\ne0$ are union of $\spec(A_{3,3,8}^{1,0})$
with two simple points, they are smoothable too. Moreover their
degree is $10$, thus they are in $\Hilb_{10}^{G,gen}(\p{4})$. When
$b=0$, the special fibre $X:={\Cal X}_0$ is defined in
$k[x_1,x_2,x_3,x_4]$ by the homogeneous ideal
$$
I:=(x_3x_4, x_2x_4, x_1x_4, x_1^2 + x_2^2, x_1x_2 + x_3^2, x_1 x_3, x_4^3 - x_1^3, x_3^3, x_2^2x_3, x_2^3).
$$
Hence it is irreducible since it is supported only on the point
$[1,0,0,0,0]\in\p4$. Moreover the corresponding algebra
$A:=k[x_1,x_2,x_3,x_4]/I\cong{\Cal A}_0$ has Hilbert function
$h_A=(1,4,4,1)$ and it is easy to check that its socle is
generated by $x_1^3$, thus $X\in{\Cal
Z}\subseteq\Hilb_{10}^{G}(\p{4})$.

We conclude that, in order to prove the irreducibility of
$\Hilb_{10}^{G}(\p{4})$, it suffices to check that
$X\not\in\Sing(\Hilb_{10}^{G}(\p{4}))$. Since
$\dim(\Hilb_{10}^{G,gen}(\p{4}))=40$ it suffices to check that the
tangent space at the point $X\in \Hilb_{10}^{G}(\p{4})$, which is
canonically identified with $H^0\big(X,{\Cal N}_{X}\big)$, has
dimension $40$.

In our case it suffices to check that
$\deg(X^{(2)})=\dim_k(k[x_1,x_2,x_3,x_4]/I^2)=50$, thanks to
Formula (2.6),  and this can be computed via any computer software
for symbolic computation. This computation concludes the proof of
the statement. \qed
\enddemo

Now assume $N\ge5$ and let $X\in {\Cal Z}_N$. Due to the
definition of ${\Cal Z}_N$ we know that there is an embedding
$X\subseteq\p4$. Thanks to the discussion above we know that $X$
is smoothable in $\p4$, thus the same holds in $\p N$ due to Lemma
2.2. This proves the following

\proclaim{Corollary 4.4} Let $X\cong\spec(A)\subseteq\p N$, where
$h_A=(1,4,4,1)$ and $N\ge 4$. Then $X$ is smoothable in $\p N$.
\qed
\endproclaim

We are now ready to summarize the results proved in this section
and in the previous one in order to give the

\demo{Proof of the Main Theorem} In order to prove that
$\Hilb_{10}^{G}(\p{N})$ is irreducible it suffices to check
$\Hilb_{10}^{G}(\p{N})=\Hilb_{10}^{G,gen}(\p{N})$, i.e. that each
Gorenstein subscheme $X\subseteq\p N$ of dimension $0$ is
smoothable.

If $X$ has at least two components this follows from Proposition
2.7. Thus we restrict our attention to irreducible schemes $X$.
Let $X\cong\spec(A)$ for some local Artinian Gorenstein
$k$--algebra with maximal ideal $\M$. If $\dim_k(\M^2/\M^3)=1$,
then the smoothability of $X$ is proven in Proposition 3.1, if
$\dim_k(\M^2/\M^3)=2$, in Propositions 3.2 and 3.3, if
$\dim_k(\M^2/\M^3)=3$ in Propositions 3.4 and 3.5, if
$\dim_k(\M^2/\M^3)=4$ in Corollary 4.4. \qed
\enddemo

Lemma 6.21 of [I--K] essentially asserts the reducibility of
$\Hilb_{14}^{G}(\p{N})$ when $N\ge6$. Indeed the authors claim the
existence of a scheme $X\cong\spec(A)\subseteq\p 6$, where
$h_A=(1,6,6,1)$ and having tangent space of dimension $76$.

Since  the main component
$\Hilb_{14}^{G,gen}(\p{N})\subseteq\Hilb_{14}^{G}(\p{N})$ has
dimension $84$ we infer the existence of a second component ${\Cal
H}\subseteq\Hilb_{14}^{G}(\p{N})$ of dimension at most $76$.

In order to construct such a scheme it suffices to make again use
of the theory of inverse systems explained above. For example, if
one considers $N=6$ and the polynomial
$$
\align
g(y_1,\dots,y_6):= y_1^3  + y_2^3 + y_3^3 + y_4^3 + y_5^3 + y_6^3 &+ (y_1+y_2+y_3+y_4+y_5+y_6)^3 +\\ +(2y_1+y_2-2y_3+y_5-y_6)^3 &+ (-y_1-2y_2-2y_3-2y_4+2y_5-2y_6)^3 +\\
&+ (-y_1-y_2+2y_3+y_4-2y_6)^3
\endalign
$$
then an explicit computation shows that the corresponding local,
Artinian, Gorenstein $k$--algebra $A$ has $h_A=(1,6,6,1)$ and,
using Formula (2.6), that
$$
h^0\big(X,{\Cal N}_{X}\big)=\dim_k(k[x_1,\dots,x_6]/(g^\perp)^2)-\dim_k(k[x_1,\dots,x_6]/g^\perp)=76.
$$

No analogous results are known for $\Hilb_{d}^{G,gen}(\p{N})$ with
$11\le d\le 13$. Similar computations with $N=5$ and polynomials
of degree $3$, give at most local, Artinian, Gorenstein
$k$--algebras $A$ with $h_A=(1,5,5,1)$ such that
$X=\spec(A)\subseteq\p6$ satisfies $h^0\big(X,{\Cal
N}_{X}\big)=60$ which is exactly the dimension of
$\Hilb_{12}^{G,gen}(\p{5})$.

For this reason we explicit the following question essentially due
to A.V. Iarrobino.

\definition{Question 4.5}
Is $\Hilb_{d}^{G}(\p{N})$  irreducible if and only if $d\le13$?
\enddefinition

\head 5. The singular locus of $\Hilb_{10}^{G}(\p{N})$
\endhead
In this last section, we describe the singular locus of
$\Hilb_{10}^{G}(\p{N})$. Since $\Hilb_{d}^{G}(\p{N})$ is
irreducible of dimension $dN$ for $d\le10$, it follows that $X$ is
obstructed, i.e. it is singular in $\Hilb_{d}^{G}(\p{N})$, if and
only if $h^0\big(\p N,{\Cal N}_X\big)>dN$.

Due to Formula (2.4) and Proposition 2.5, this can happen only
when there is an irreducible component $Y\subseteq X$ of degree
$d$ in the following list: \roster \item $Y\cong\spec(A_{n,d})$,
with $6\le n+2\le d$; \item $Y\cong\spec(A_{n,2,d}^t)$ with
$t=1,2$ and $8\le n+4\le d$; \item
$Y\cong\spec(A_{n,3,n+5}^{t,\alpha})$ with $t=1,\dots,6$ and $9\le
n+5=d$  (if $n=5$ then $Y=X$); \item
$Y=X\cong\spec(A_{4,3,10}^{t})$ with $t=0,\dots,6$; \item
$Y=X\cong\spec(A_{4,2,2,10}^{t})$ with $t=1,2,3$; \item
$Y=X\in{\Cal Z}_N$.
\endroster

In Section 5 of [C--N2] we checked that in cases (1) and (2), the
corresponding schemes are obstructed. In case (3) it is proven
there that $Y$ is obstructed if and only if $t=4,5,6$ when $n=4$,
by computing explicitly $h^0\big(Y,{\Cal N}_Y\big)$ where the
embedding $Y\subseteq\a n\subseteq \p n$ is the natural one
corresponding to the representation of $Y$ as spectrum of a
quotient of $k[x_1,\dots,x_n]$ and making use of Formula (2.6) as
already done above. We now examine with the same approach, using
any computer software for symbolic calculations, the cases (3)
with $n=5$ and (4), (5) with $n=4$.

In case (3) we have that the normal sheaf ${\Cal N}_X$ of the
embedding induced by the natural quotient
$k[x_1,\dots,x_5]\twoheadrightarrow A_{n,3,n+5}^{t,\alpha}$
satisfies
$$
h^0\big(X,{\Cal N}_{X}\big)=\cases
57\quad&\text{if $t=1,2,3$},\\
64\quad&\text{if $t=4,5,6$}.
\endcases
$$
In case (4) with respect to the natural quotient
$k[x_1,\dots,x_4]\twoheadrightarrow A_{4,3,10}^{t}$ we have
$$
h^0\big(X,{\Cal N}_{X}\big)=\cases
40\quad&\text{if $t=0,1$},\\
45\quad&\text{if $t=2,3,4,5,6$}.
\endcases
$$
Finally, in case (5), with respect to the natural quotient
$k[x_1,\dots,x_4]\twoheadrightarrow A_{4,2,2,10}^{t}$ we have
$$
h^0\big(X,{\Cal N}_{X}\big)=45\quad\text{if $t=1,2,3$}.
$$

We can summarize the above results in the following

\proclaim{Theorem 5.3} Let $X\in
\Hilb_{10}^{G}(\p{N})\setminus{\Cal Z}_N$. Then $X$ is obstructed
if and only if it contains an irreducible component isomorphic to
either $\spec(A_{n,d})$ or $\spec(A_{n,2,d}^t)$, where $n\ge4$, or
$\spec(A_{4,3,9}^{t,\alpha})$, where $t=4,5,6$, or
$\spec(A_{4,3,10}^{t})$, where $t=2,3,4,5,6$, or
$\spec(A_{5,3,10}^{t,\alpha})$  or $\spec(A_{4,2,2,10}^{t})$,
without restrictions on $t$. \qed
\endproclaim
It is natural to ask what happens in the case $X\in{\Cal
Z}_N\subseteq\Hilb_{10}^{G}(\p{N})$. We checked in the previous
section that the general scheme in ${\Cal Z}_N$ is not obstructed.
In principle the theory of inverse system and the classification
of cubic surfaces (e.g. as the one in [B--L]) could allow us to
complete the description of points in ${\Cal Z}_N$, hence it could
help to describe completely the singular locus of
$\Hilb_{10}^{G}(\p{N})$.

Unfortunately, taking into account of [B--L], we have at least
$22$ different cases to handle, most of them depending on many
parameters. Thus a direct approach seems to be useless in this
case. Thus we need another method. Notice that each point in
${\Cal Z}_N$ corresponds to a local Artinian Gorenstein
$k$--algebra $A$ with Hilbert function $(1,4,4,1)$.

As explained in the previous section, such kind of algebra is
naturally graded, i.e. it can be written as a suitable quotient
quotient $S:=k[x_1,x_2,x_3,x_4]/I$ with $I$ homogeneous and it
corresponds, via Macaulay's correspondence, to a cubic form $g$,
i.e. $I=g^\perp$.

\proclaim{Lemma 5.5} The minimal free resolution of $A\cong
S/g^\perp$ over $S$ has the form
$$
\align
0  \longrightarrow S(-7)  \longrightarrow S^6(-5) \oplus S^\beta(-4)  \longrightarrow& S^{5+\beta}(-4) \oplus S^{5+\beta}(-3) \longrightarrow \\
\longrightarrow& S^\beta(-3) \oplus S^6(-2)  \longrightarrow S  \longrightarrow  A  \longrightarrow  0
\endalign
$$
for some $\beta\ge0$.
\endproclaim
\demo{Proof} The ideal $g^\perp$ has obviously six minimal
generators of degree $2$, but it could also have some more minimal
generators in degree $3$ or higher. Thus  the  minimal free
resolution of $A$ over $S$ ends with
$$
S^6(-2)\oplus S(-3)^{\beta}  \oplus F\longrightarrow S  \longrightarrow A  \longrightarrow 0
$$
where $\beta\ge0$ is the number of minimal cubic generators of
$g^\perp$ and $F$ is a direct sum of $S(-j)$ with $j\ge4$.

If $F$ does not contain the direct summand $S(-4)$,then the cubic
forms in the ideal $g^\perp$ would generate its degree $4$
homogeneous part, thus they would generate $g^\perp$ in degree
greater than $3$, i.e. $F=0$. It remains to examine the case when
$g^{\perp} $ has a minimal generator in degree $ 4$.

Since $A$ is Gorenstein with maximum socle degree $3$, it follows that the
minimal free resolution is self--dual up to twisting by $S(-7)$
(this is a well--known fact. For the sake of completeness we quote
[B--H] as reference: in particular Corollary 3.3.9, Proposition
3.6.11,  Examples 3.6.15, Theorem 3.6.19 and the remark after it).
Moreover the middle free module cannot contain $S(-2)$ summands
since the generators in degree $2$ are obviously linearly
independent. Combining such remarks we obtain that the minimal
free resolution of $A$ has the shape
$$
\align
0  \longrightarrow S(-7)  \longrightarrow& S^6(-5) \oplus S^\beta(-4) \oplus \check F(-7)\longrightarrow G\longrightarrow \\
\longrightarrow& S^\beta(-3) \oplus S^6(-2)  \oplus F\longrightarrow S  \longrightarrow  A  \longrightarrow  0
\endalign
$$
On one hand, by assumption, $S(-4) $ is a free addendum of $ F$,
hence $ S(-3) $ is a free addendum of $\check F(-7)$. On the other
hand  the resolution above is minimal, thus at each step the
minimal degree of the syzygies must increase at least by one. This
two remarks yields a contradiction, thus $F=0$.

For the same reasons $G$ contains only direct summand of the form
$S(-j)$ with $j\ge3$ and $G\cong \check G(-7)$. A simple
computation thus yields $G\cong S^{5+\beta}(-4) \oplus
S^{5+\beta}(-3)$. \qed
\enddemo

\remark{Remark 5.6} Notice that the argument above can be also
used for proving the following assertion. Let $I\subseteq
k[x_1,\dots,x_N]$ be a homogeneous ideal such that
$A:=k[x_1,\dots,x_N]/I$ is a local Artinian Gorenstein
$k$--algebra with maximum socle degree $e$. Then $I$ has a minimal generator in
degree $e+1$ if and only if $A\cong k[t]/(t^{e+1})$ or,
equivalently, if and only if $I=g^\perp$ with $g=\ell^{e+1}$ for
some linear form $\ell\in k[y_1,\dots,y_N]$.
\endremark
\medbreak

At this point we are ready to start with our classifications
results. We first examine the general case.

\proclaim{Proposition 5.7}
Using the notation above let $ A^{(2)}:=S/(g^{\perp})^2 $.
If $\beta=0$ in Lemma 5.5, then $ h_{A^{(2)}}=(1,4,10,20,14,1)$.
\endproclaim
\demo{Proof} Let $ f_1, \dots, f_6 \in S_2$ be a minimal set of
quadratic generators of $g^\perp$. Since the ring $ A $ is
Artinian, we can assume that $ f_1,\dots, f_4 $ is a regular
sequence in $ S$. To fix notation, we assume that the first map $
\varphi: S^6(-2) \to S $ of the resolution in Lemma 5.5 of $ A $
is given by $ \varphi(e_i) = f_i $ for $ i=1, \dots, 6,$ where $
e_1, \dots, e_6 $ is the canonical basis of $S^6(-2)$.

Let $ M := (M_1 \vert M_2) $ be the matrix representing the map $
S^5(-4) \oplus S^5(-3) \to S^6(-2) $ with respect to the canonical
bases of the involved free modules. Trivially the elements of $
M_1 $ have degree $2$ while the ones of $ M_2 $ have degree $ 1$.
Let $ V \subseteq S_1 $ (resp. $ W \subseteq S_1$) be the subspace
generated by the elements of the $ 5$--th row (resp. $ 6$--th row)
of $ M_2$. If $ \dim_k(V) \le 2 $ and $ \dim_k(W) \le 2$, then we
can obtain a degree $ 1 $ syzygy of $ g^\perp$ with the last two
entries equal to $ 0$, that is to say, there exists a degree $ 1 $
syzygy of $ f_1, \dots, f_4$, a contradiction, since the
resolution of  $ I =(f_1, \dots, f_4)\subseteq S $ is Koszul being
$ f_1, \dots, f_4$ a regular sequence. Hence, either $ V $ or $ W
$ has dimension at least $3$. Up to exchange $ f_5 $ and $ f_6$,
we can finally assume $ \dim_k(W) \ge 3$.

The minimal free resolution of $ B := S/I $ is
$$
0\longrightarrow  S(-8)\longrightarrow  S^4(-6) \longrightarrow  S^6(-4) \longrightarrow S^4(-2) \longrightarrow  S\longrightarrow  B \longrightarrow  0,
$$
whence $ h_B = (1,4,6,4,1)$. Since $I\subseteq g^\perp$, it
follows the existence of a natural epimorphism $ B
\twoheadrightarrow A$ with kernel $g^\perp/I$.

Of course, the classes of $f_5 $ and $ f_6 $ mod $ I $ are in $
B_2$. It is then obvious that $ f_5 S_d \subset I$ and $ f_6 S_e
\subset I $ for some integers $ d, e$. Let $ J :=( f_1, \dots,
f_5)$.

We first consider the case $ \dim_k(W) = 4$, i.e. $ W = S_1$. In
this case $ f_6 S_1 \subset J$. From the above inclusion and the
short exact sequence
$$
0 \longrightarrow g^\perp/ J\longrightarrow S/ J \longrightarrow A\longrightarrow 0
$$
we deduce $ h_{S/ J} = (1,4,5,1)$. Consider now the exact sequence
$$
0  \longrightarrow  J/ I  \longrightarrow S /I  \longrightarrow S/J \longrightarrow 0.
$$
By computing the dimensions of the homogeneous pieces, we obtain
$ \dim_k( (J/I)_j )=1, 3, 1,$ for $ j = 2, 3, 4,$ respectively,
and $ 0 $ otherwise. Hence, there exists $ \ell_1 \in S_1 $ such
that $ \ell_1 f_5 \in I$, and, if $ \ell_1, \dots, \ell_4 $ is a
basis of $ S_1$, we infer that the cosets of $ \ell_2 f_5, \ell_3
f_5, \ell_4 f_5 $ are linearly independent in $ S / I$.

Looking at the matrix $ M_2$, after reducing its columns by
elementary operations, we can say that there is one column whose
last two entries are $ \ell_1,0,$ respectively. After reducing the
columns of $ M $ by elementary operations, all the elements of the
$ 5$--th row of $M_1 $ are non-zero. Hence, there are $ 5 $
linearly independent elements in $ (\ell_2, \ell_3, \ell_4)^2 $
which are in $ I $ when multiplied by $ f_5$. Since there are no
minimal syzygies in degree $ 3 $ and $ (\ell_2, \ell_3, \ell_4)^3
f_5\subseteq I$, we can choose generators of $ (\ell_2, \ell_3,
\ell_4) $ in such a way that $ \ell_2^2 f_5, \ell_2 \ell_3 f_5,
\ell_3^2 f_5, \ell_2 \ell_4 f_5, \ell_4^2 f_5 \in I$, while the
coset of $\ell_3 \ell_4 f_5 $ spans $ J/I $ in degree $ 4$.

The ideals $ I,J, g^\perp $ give rise to the following sequence of
strict inclusions
$$
I^2 \subset IJ \subset J^2 \subset Jg^\perp\subset (g^\perp)^2
$$
that we will use in order to compute $h_{A^{(2)}}$.

To start with, we consider $ B^{(2)} = S / I^2$. On one hand it
fits into the exact sequence
$$
0 \longrightarrow  I/I^2\longrightarrow   B^{(2)} \longrightarrow   B \longrightarrow   0.
$$
On the other hand $ I/I^2= I \otimes_SS/I \cong (S/I)^4(-2)$,
since $I$ is generated by a regular sequence of quadratic forms.
Hence $h_{B^{(2)}} = (1,4,10,20,25,16,4)$.

The module $ IJ / I^2 $ is generated by the cosets of $ f_1 f_5,
\dots, f_4 f_5$. Let $ a_1, \dots, a_4 \in S $ be such that $
a_1f_1 f_5 + \dots + a_4 f_4 f_5 \in I^2.$ Hence, $ (a_1 f_5,
\dots, a_4 f_5) $ is zero in $ (S/I)^4(-2)$, i.e. $ a_i f_5 \in I
$ for each $ i = 1,\dots, 4$. Thanks to the previous discussion,
this happens if, and only if $ a_1 = 0 $ (when $ \deg (a_i) = 0$),
$ a_i \in (\ell_1) $ (when $\deg (a_i) = 1$), $ a_i \in (\ell_1,
\ell_2^2, \ell_2 \ell_3,\ell_3^2, \ell_2 \ell_4, \ell_4^2) $ (when
$ \deg (a_i )= 2$) and, finally, $ a_i \in S_j $ (when $ \deg (a_i
)= j \ge 3$). Hence  $ \dim_k ((IJ / I^2)_j) = 4, 12, 4$, for $ j
= 4, 5, 6$ respectively, and $ 0 $ otherwise, thus $ h_{S / IJ} =
(1,4,10,20,21,4)$.

Now, consider $C^{(2)}:=S/J^2$ and the exact sequence
$$
0 \longrightarrow J^2/IJ\longrightarrow  S/IJ \longrightarrow  C^{(2)} \longrightarrow 0.
$$
The module $J^2 / IJ$ is generated by the coset of $ f_5^2,$ and
the assertion $ a f_5^2\in IJ$ is equivalent to the assertion $ a
f_5 \in I$. It thus follows from the above discussion and from the
computation of $h_{S / IJ}$ we get that $\dim_k( (J^2 / IJ)_j)= 1,
3$, for $ j= 4, 5,$ respectively, and $ 0 $ otherwise. Hence $
h_{C^{(2)}} =(1,4,10,20,20,1)$.

The module $ Jg^\perp / J^2 $ is generated by the cosets of $ f_1
f_6, \dots, f_5 f_6$, thus the dimensions of its homogeneous
pieces are $ \dim_k ((Jg^\perp / J^2)_j )= 5 $ if $ j = 4$, and $
0 $ otherwise, since $ f_6 S_1 \subseteq J$. Hence the the Hilbert
function of $S/Jg^\perp$ can be computed by using the exact
sequence
$$
0  \longrightarrow  Jg^\perp/J^2 \longrightarrow C^{(2)}
\longrightarrow  S/Jg^\perp \to 0.
$$
We obtain $h_{S /Jg^\perp}= (1,4,10,20,15,1) $.

Finally, the module $ (g^\perp)^2/Jg^\perp$ is generated by the
coset of $f_6^2 $ and it is non-zero only in degree $ 4$. The
Hilbert function of $ A^{(2)} $ is then equal to $ h_{A^{(2)}} =
(1,4,10,20,14,1)$ as it comes from considering the exact sequence
$$
0 \longrightarrow (g^\perp)^2/Jg^\perp \longrightarrow  S/Jg^\perp \longrightarrow  A^{(2)}  \longrightarrow  0.
$$
Thus the statement is proved under the extra hypothesis
$\dim_k(W)=4$.

Now, we consider the case $ \dim_k(W)= 3$. Of course, up to
exchanging the roles of $f_5$ and $f_6$, we can also assume $
\dim_k(V)\le 3$. We can reduce the matrix $ M_2 $ by using
elementary operations on its columns, and so we can assume that
two entries of the $ 6$--th row of $ M_2 $ are equal to $ 0$.
Moreover, from the three non--zero entries of the row, we deduce
that $ \ell f_6 \in J $ for each $ \ell\in W$ and that the last
two columns of $ M_2 $ have two linearly independent elements on
the $ 5$--th row.

Recall that $\dim_k(V)$ is either $2$ or $3$. In the former case
we can assume that, for each column of $ M_2$, if the element of
the $ 5$--th row is non--zero, then the element on the $ 6$--th
row is zero and conversely. In the latter case we can assume that
the previous situation happens on $ 4 $ columns of $ M_2$.
Furthermore, if we reduce the matrix $ M $ by using elementary
operations on its columns, not all the entries of the $ 6$--th row
of $ M_1 $ can be equal to $0$, due to the fact that $ f_6 S_e
\subseteq I$.

Then, if $ S_1 = W \oplus\langle \ell \rangle$, we can assume that
$ \ell^2 f_6 \in J$. Let $C =S/J$ and consider the short exact
sequence
$$
0 \longrightarrow g^\perp/J \longrightarrow S/J  \longrightarrow A  \longrightarrow0
$$
From the discussion above, we deduce that $ h_{S/J}(1) = 4$ and $
h_{S/J}(2) = 5$. Moreover, $ \ell f_6 \notin J $ and so $
h_{S/J}(3) = 2$, but $ h_{S/J}(4) = 0$, since  $ \ell^2 f_6 \in
J$. Hence $ h_{S/J} = (1,4,5,2)$. We can also consider the short
exact sequence
$$
0\longrightarrow J/I\longrightarrow S/I\longrightarrow S/J\longrightarrow 0.
$$
Thus the Hilbert function of $ J/I $ satisfies $\dim_k( (J/I)_j )=
1, 2, 1 $ for $ j = 2, 3, 4$ respectively, and $ \dim_k( (J/I)_j
)= 0 $ otherwise.

From the analysis of the elements of the $ 5$--th row of $ M $
corresponding to the  $ 0 $ entries on the last row of $ M$, we
get that there exists a dimension $ 2 $ subspace $ V' \subseteq
S_1 $ such that $ \ell f_5 \in I $ for each $ \ell\in V'$. Let us
choose $ V'' \subseteq S_1 $ such that $ S_1 = V' \oplus V''$. Let
$ \ell_1, \ell_2 $ be a basis of $ V''$. Then $ J/I $ is generated
by the coset of $ f_5 $ in degree $ 2$ and by the cosets of $
\ell_i f_5, i = 1,2$, in degree $ 3$. Furthermore, we have that
two among $\ell_1^2 f_5, \ell_1 \ell_2 f_5, \ell_2^2 f_5 $ are in
$ I$. The columns of $ M_2 $ have degree $ 2$, and so $ \ell_i^2
f_5 \in I$, $i = 1, 2$, since $ f_5 (\ell_1, \ell_2)^3 \subseteq
I$, but we have no minimal syzygies in degree $ 3$.

As in the case $ \dim_k(W)= 4$, the ideals $ I$, $J$ and $g^\perp$
give rise to the following sequence of strict inclusions
$$
I^2 \subset IJ \subset J^2 \subset Jg^\perp\subset (g^\perp)^2
$$
that we will use again to compute $h_{A^{(2)}}$. The Hilbert
function of $ B^{(2)} $ has been already computed above, and we do
not repeat the computation.

The module $ IJ / I^2 $ is generated by the cosets of $ f_1
f_5,f_2, f_5, f_3 f_5, f_4 f_5 $ and fits into the short exact
sequence
$$
0 \longrightarrow  IJ/I^2\longrightarrow  B^{(2)} \longrightarrow S/IJ\longrightarrow 0.
$$
Let $ a_1, \dots, a_4 \in S $ be such that $ \sum_{i=1}^4 a_if_i
f_5 \in I^2$. Then $ (a_1 f_5, \dots, a_4 f_5) $ is zero in $
(S/I)^4(-2)$, i.e. $a_i f_5 \in I$ for each $ i = 1, \dots, 4$. If
$ \deg (a_i) = 0$, this implies $ a_i = 0$, for each $ i$; if $
\deg(a_i) = 1$, we get $ a_i \in V' $ for each $ i$; if $ \deg
(a_i) = 2$, then $ a_i \in V' S_1 + \langle \ell_1^2, \ell_2^2
\rangle$; finally, if $ \deg(a_i) \ge 3$, then $ a_i f_5 \in I $
for each $ i$. It follows that $ \dim_k((IJ/I^2)_j)= 4, 8, 4$, for
$ j = 4, 5, 6$ respectively and  $ 0 $ otherwise. Hence $ h_{S
/IJ} = (1,4,10,20,21,8)$.

The next step consists in considering the short exact sequence
$$
0\longrightarrow J^2/IJ\longrightarrow S/IJ\longrightarrow  C^{(2)} \longrightarrow  0.
$$
The module $ J^2 / IJ $ is generated by the coset of $ f_5^2$. We
know that  $ S_j = (IJ)_j $ for $ j \ge 6$, hence it is enough to
consider $ a \in S $ such that $ a f_5^2 \in IJ$, with $
\deg(a)\le 1.$ This means that $ a f_5 \in I$, and so either $ a=
0 $ (when $ \deg(a)= 0$) or $ a \in V' $(when $ \deg(a)= 1$). It
follows that $ \dim_k((J^2/IJ)_j) = 1, 2$, for $ j = 4, 5$
respectively, and $ 0 $ otherwise. Hence, the Hilbert function of
$ C^{(2)} $ is $ h_{C^{(2)}} = (1,4,10,20,20,6)$.

The module $ IJ/J^2 $ is generated by the cosets of $ f_1 f_6,
\dots, f_5 f_6$. Then, we have that $ \dim_k((Jg^\perp/J^2)_4 )=
5$. Let $ a \in S_1$, and consider $ a f_5 f_6$. If $ a \in W$,
then $ a f_6 \in J$, and so $ a f_5 f_6 \in Jg^\perp$. If $ a \in
V'$, then $ a f_5 \in I$, and so $ a f_5 f_6 \in (f_1 f_6, \dots,
f_4 f_6)$. Hence, if $ W + V' = S_1$, we deduce that $
(Jg^\perp/J^2)_5$ is spanned by the cosets of $ \ell f_1 f_6,
\dots, \ell f_4 f_6 $, whence $ \dim_k((Jg^\perp/J^2)_5) =4$. If $
W \supset V'$, then the cosets of $ \ell f_1 f_6, \dots, \ell f_5
f_6 $ are linearly independent, thus $ \dim_k((Jg^\perp/J^2)_5)=
5$. Hence, the Hilbert function of $ S / Jg^\perp $ is either $
h_{S / Jg^\perp } = (1,4,10,20,15,2) $ (when $ V'\not\subseteq
S_1$) or $h_{S / Jg^\perp } =(1,4,10,20,15,1)$, (when $ V' \subset
W$), as we easily obtain  from the short exact sequence
$$
0 \longrightarrow Jg^\perp/J^2\longrightarrow C^{(2)} \longrightarrow S/Jg^\perp\longrightarrow 0.
$$
In both the cases, $ (g^\perp)^2 / Jg^\perp $ is generated by the
coset of $ f_6^2$, hence $ \dim_k(((g^\perp)^2/Jg^\perp)_4) = 1$.

If $ V' \not\subseteq W$, then $\ell f_6^2 $ spans $
((g^\perp)^2/Jg^\perp)_5 $ as vector space, thus
$\dim_k(((g^\perp)^2/Jg^\perp)_5) = 1$. From the exact sequence
$$
0  \longrightarrow (g^\perp)^2/Jg^\perp\longrightarrow S/Jg^\perp  \longrightarrow  A^{(2)}  \longrightarrow  0
$$
we finally obtain that $ h_{A^{(2)}} = (1,4,10,20,14,1)$. If $ V'
\subset W,$ then we certainly have $ h_{A^{(2)}} =
(1,4,10,20,15,1)-(0,0,0,0,1,h_5)=(1,4,10,20,14,1-h_5)$ where
$h_5\ge0$. Due to the Main Theorem then $\Hilb^G_{10}(\p4)$ is
irreducible, thus the scheme $X:=\spec(A)$ embedded in
$\a4\subseteq\p4$ via the natural quotient $S\twoheadrightarrow A$
lies in a scheme of dimension $40$. Proposition 2.5 thus yields
that
$$
40\le h^0\big(X,{\Cal N}_X\big)=\dim_k(A^{(2)})-\dim_k(A)=40-h_5,
$$
whence $h_5=0$. We conclude that $ h_{A^{(2)}} = (1,4,10,20,14,1)$
also in this second case. \qed
\enddemo

Now we examine the case when $\beta\ge1$. In this case $g^\perp$
has at least one minimal cubic generator. By [C-R-V], Theorem
6.18, there exists $ \ell \in S_1 $ such that $ \ell(g) \in R_2 $
is a rank $ 1 $ quadric. Up to a change of coordinates, we can
assume $ \ell = x_4$, and $ x_4 (f )= y^2 $ for some $ y =
\sum_{i=1}^4 b_i y_i \in R_1$.  Either $b_4 \ne 0 $ or $b_4 = 0$.

In the former case we can assume that $b_4=1$. If $ b_1 y_1 + b_2
y_2 +b_3 y_3 = 0$, then $ g=y_4^3 + g_0 $ for a suitable $ g_0 \in
k[y_1, y_2, y_3]$. If $ b_1 y_1 + b_2 y_2 +b_3 y_3 $ is non--zero,
then, up to a change of variables, we have $g =y_4^3 + y_4^2 y_2 +
y_4 y_2^2 + {g}_1 $ for a suitable cubic form $ g_1 \in k[y_1,
y_2, y_3]$. By setting $ x_2 = X_4-X_2, x_i = X_i$ for $ i=1, 3,
4,$ and $ y_4 = Y_4 + Y_2$, $y_2 = -Y_2$, $y_i = Y_i $ for $ i =
3, 4$, then $ g = Y_4^3 + Y_2^3 +g_2(-Y_1, Y_2, Y_3) $.

In the latter case we have that $b_4=0$. Necessarily $ b_1 y_1 +
b_2 y_2 +b_3 y_3 \ne0$, hence up to a proper change of the
variables we can assume  $g= y_3^2y_4 + \widehat{g}(y_1, y_2,
y_3)$.

The above discussion proves the \lq\lq only if\rq\rq\ of the
following

\proclaim{Lemma 5.8} 
Let $g\in S_3$. Then, $ g^\perp $ has minimal
generators in degree $ 3 $ if and only if there exists a cubic
form $ \widehat{g} \in k[y_1, y_2, y_3]$ such that, up to a proper
choice of coordinates in $ R$, either $ g=y_4^3 + \widehat{g} $ or
$ g = y_3^2y_4 + \widehat{g}$.
\endproclaim
\demo{Proof} It remains to prove the \lq\lq if\rq\rq\ part. If
either $ g=y_4^3 + \widehat{g} $ or $ g = y_3^2y_4 + \widehat{g}$
for some cubic form $\widehat{g} \in k[y_1,y_2, y_3]$, then $ x_4
(g)$ is equal either to $ 3y_4^2 $ or to $ 2y_3^2$, hence $
x_4(g)$ is a rank $ 1 $ quadric. Again by [C-R-V], Theorem 6.18,
$g^\perp$ has a minimal generator in degree $ 3$. \qed
\enddemo

We now go to complete our classification.

\proclaim{Proposition 5.9} 
Using the notation above let $
A^{(2)}:=S/(g^{\perp})^2 $. If $\beta\ge1$ in Lemma 5.5, then
either $\beta=1$ and $ h_{A^{(2)}}=(1,4,10,20,14,1)$ or $\beta=3$
and $ h_{A^{(2)}}=(1,4,10,20,16,4)$.
\endproclaim
\demo{Proof} 
Due to Lemma 5.8 we can assume that either $ g=y_4^3
+ \widehat{g} $ or $ g = y_3^2y_4 + \widehat{g}$ for some cubic
form $\widehat{g} \in k[y_1,y_2, y_2]$.

Consider the first case. Up to a suitable change of coordinates,
$\widehat{g} $ is equal to one of the following
$$
\gather
y_3^3,\qquad
y_2y_3^2,\qquad y_2y_3(y_2-y_3),\qquad y_1y_2y_3, \qquad y_3(y_1y_3-y_2^2),\qquad
y_2(y_1y_3-y_2^2),\\
y_1^2y_3+y_2^2y_3-y_2^3, \qquad y_1^2y_3-y_2^3,\qquad
y_1^2y_3-y_2^3+(1+t)y_2^2y_3-ty_2y_3^2
\endgather
$$
where $ t \in k $ is
different from $ 0 $ and $ 1$. In the
various cases we perform the computation using any computer software for symbolic calculations, and we report the results.

The first three choices give Artinian Gorenstein rings with
Hilbert function different from $ (1,4,4,1),$ because $ g $ is a
cone in those cases.

If $ g = y_4^3 + y_1y_2y_3,$ then
$$
g^\perp = (x_1^2, x_2^2,
x_3^2, x_1x_4, x_2x_4, x_3x_4, 6x_1x_2x_3-x_4^3).
$$
Hence, $\beta=1$, and $ h_{A^{(2)}} = (1,4,10,20,14,1)$.

If $ g = y_4^3 + y_3(y_1y_3-y_2^2)$, then
$$
g^\perp = (x_1^2,
x_1x_2, x_2^2+x_1x_3, x_1x_4, x_2x_4, x_3x_4, 3x_1x_3^2-x_4^3,
x_2x_3^2, x_3^3).
$$
Hence, $ \beta=3$, and $ h_{A^{(2)}} =
(1,4,10,20,16,4)$.

If $ g = y_4^3 + y_2(y_1y_3-y_2^2),$ then
$$
g^\perp = (x_1^2,
x_2^2+6x_1x_3, x_3^2, x_1x_4, x_2x_4, x_3x_4, 6x_1x_2x_3-x_4^3).
$$
Hence, $\beta=1$, and $ h_{A^{(2)}} = (1,4,10,20,14,1)$.

If $ g = y_4^3 + y_1^2y_3+y_2^2y_3-y_2^3,$ then
$$
g^\perp =
(x_1^2-x_2^2-3x_2x_3, x_1x_2, x_3^2, x_1x_4, x_2x_4, x_3x_4,
3x_2^2x_3-x_4^3).
$$
Hence, $\beta=1,$ and $ h_{A^{(2)}} =
(1,4,10,20,14,1).$

If $ g = y_4^3 + y_1^2y_3-y_2^3,$ then
$$
g^\perp = (x_1x_2,
x_2x_3, x_3^2, x_1x_4, x_2x_4, x_3x_4, x_1^3, x_2^3+x_4^3,
3x_1^2x_3-x_4^3).
$$
Hence, $ \beta=3,$ and $ h_{A^{(2)}} =
(1,4,10,20,16,4).$

If $ g = y_4^3 + y_1^2y_3-y_2^3+(1+t)y_2^2y_3-ty_2y_3^2,$ then
$$
\align
g^\perp = (x_1x_2, x_1x_4, &x_2x_4, x_3x_4,
t(1+t)x_1^2-tx_2^2+3x_3^2, \\
(t^2-t&+1)x_1^2-(1+t)x_2^2-3x_2x_3,x_1^3, x_1x_3^2, x_3^3, x_2^3+x_4^3, \\
&tx_1^2x_3+x_2x_3^2,
(1+t)x_1^2x_3-x_2^2x_3, 3x_1^2x_3+x_2^3).
\endalign
$$
If $ t^2-t+1 \not= 0,$
then $ \beta=1,$ and $ h_{A^{(2)}} = (1,4,10,20,14,1).$ If $
t^2-t+1=0,$ then $ \beta=3 $ and $ h_{A^{(2)}} =
(1,4,10,20,16,4).$ By the way, it is well--known that the condition $ t^2-t+1=0 $ corresponds to the $
j$--invariant of the smooth cubic to be $ 0$ (see [Ha], Section IV.4).

Let us consider now the second case, i.e. $ g = y_3^2y_4 +
\widehat{g}$ for some cubic form $\widehat{g} \in k[y_1,y_2,
y_2]$. With a change of coordinates, we can assume that
$$
g =
y_3^2y_4 + y_3(b_1y_2^2+2b_2y_1y_2+b_3y_1^2) +
(b_4y_2^3+b_5y_1y_2^2+b_6y_1^2y_2+b_7y_1^3).
$$

The form $b_4y_2^3+b_5y_1y_2^2+b_6y_1^2y_2+b_7y_1^3$ in the expression of $g$ can have either three simple roots, or a triple root, or a simple root and a double one. According to its roots,  up to a change of coordinates, it can be written as either $ y_1^3 + y_2^3$, or $ y_2^3$, or $ y_1y_2^2$. Accordingly $g$ has one of the following forms:
$$
\gather
y_3^2 y_4 + y_3(b_1y_2^2 + 2b_2y_1y_2+ b_3y_1^2) + y_1^3+y_2^3,\qquad  y_3^2 y_4 + y_3(b_1y_2^2+2b_2y_1y_2+ b_3y_1^2) + y_2^3,\\
y_3^2 y_4 + y_3(2b_2y_1y_2+ b_3y_1^2) + y_1y_2^2
 \endgather
 $$
(in the last case we made the extra change of variables $y_\mapsto y_1+b_1y_3$).

In the first case, we have that 
$$
\align
g^\perp = (x_4^2, &x_2x_4,
x_1x_4, x_1x_2-b_2x_3x_4,
3x_2x_3-b_2x_1^2-b_1x_2^2+(b_1^2+b_2b_3)x_3x_4,\\
&3x_1x_3-b_3x_1^2-b_2x_2^2+(b_1b_2+b_3^2)x_3x_4, x_1^3-3x_3^2x_4,
x_2^3-3x_3^2x_4, x_3^3).
\endalign
$$
If $ b_2 \not= 0,$ then $ \beta = 1 $
and $ h_{A^{(2)}} = (1,4,10,20,14,1).$ If $ b_2 = 0,$ then $ \beta
= 3,$ and $ h_{A^{(2)}} = (1,4,10,20,16,4).$

In the second case, we have that
$$
\align
g^\perp = (x_4^2, x_2x_4,&x_1x_4, x_1^2-b_3x_3x_4, x_1x_2-b_2x_3x_4, \\
3b_2x_1x_3 &- 3b_3x_2x_3 + (b_1b_3-b_2^2)x_2^2 - b_1(b_1b_3-b_2^2)x_3x_4, \\
&x_3^3, x_1x_3^2, x_2x_3^2, x_2^3-3x_3^2x_4, x_2^2x_3-b_1x_3^2x_4).
\endalign
$$
If $ b_2 = b_3 = 0,$ then $ g $ is a cone, and so $ g^\perp $ is degenerate. Hence, we can assume that either $ b_2 \not= 0,$ or $ b_3 \not= 0.$ In both cases, $ \beta = 3,$ and $ h_{A^{(2)}} = (1,4,10,20,16,4).$

In the last case, we have that
$$
\align
g^\perp = (x_4^2, &x_2x_4, x_1x_4, x_1^2-b_3x_3x_4, x_2x_3-b_2x_2^2, \\
&x_1x_3-b_2x_1x_2-b_3x_2^2+b_2^2x_3x_4, x_3^3, x_2^3, x_1x_3^2, x_1x_2^2-x_3^2x_4).
\endalign
$$
If $ b_3 \not= 0,$ then $ \beta=1 $ and $ h_{A^{(2)}} = (1,4,10,20,14,1)$. If $ b_3 = 0,$ then $ \beta=3 $ and $ h_{A^{(2)}} = (1,4,10,20,16,4)$.
\qed
\enddemo

Now let $ g \in R_3 $ and $A:=S/g^\perp$. Let
$X:=\spec(A)\subseteq\a4\subseteq\p4$ be the embedding associated
to the quotient $k[x_1,x_2,x_3,x_4]\twoheadrightarrow A$. An
immediate consequence of Propositions 5.6 and 5.7, of Formula
(2.4) and of Proposition 2.5 is that the normal bundle ${\Cal
N}_X$ satisfies
$$
h^0\big(X,{\Cal N}_X\big)=\cases
40&\text{if $g$ is as in Proposition 5.7,}\\
45&\text{if $g$ is as in Proposition 5.8.}
\endcases
$$
The same argument used in the proof of Theorem 5.3, thus yields

\proclaim{Theorem 5.9} Let $ g \in R_3 $, $A:=S/g^\perp$ and
$X:=\spec(A)\in{\Cal Z}_N\subseteq \Hilb_{10}^{G}(\p{N})$.  The
scheme $X$ is obstructed if and only if $\beta = 3$. \qed
\endproclaim

\enddocument